\documentclass{article}

\newtheorem{fed}{\textbf{Definition}}[section]
\newtheorem{thm}[fed]{\textbf{Theorem}}
\newtheorem{lemma}[fed]{\textbf{Lemma}}

\newtheorem{cor}[fed]{\textbf{Corollary}}
\usepackage{amssymb,bbm,graphicx,epsfig,psfrag,epic,eepic,latexsym}
\usepackage{amsmath}
\usepackage{mathrsfs}
\usepackage[dvips]{color}

\begin{document}
\title{On Kepler's geometric approach to consonances}
\author{Urs Frauenfelder}
\maketitle

\begin{abstract}
Kepler's thinking is highly original and the inspiration for discovering his famous third law is based on his rather curious geometric approach in his Harmonices mundi for explaining consonances. In this article we try to use a modern mathematical approach based on Kepler's ideas how to characterize the seven consonances with the help of the numbers of edges of polygons constructible by ruler and compass.   
\end{abstract}

\section{Introduction}

Kepler's "World harmony" (Harmonices mundi) \cite{kepler2} from 1618 is famous, because in this book Kepler states his third law on planetary motion, namely that the cubes of the semi-major axis of the elliptic orbits of the planets are proportional to the squares of their periods. However, this crucial discovery cannot found until the fifth book and in fact it remains rather mysterious how Kepler actually found his third law. 
\\ \\
Kepler's thinking is highly original and rather unique. One characteristic feature of Kepler's thinking is that he is looking for geometrical and not arithmetical explanations. Kepler's geometric approach played already a crucial role in his "Cosmographic Mystery" (Mysterium Cosmographicum) \cite{kepler0} from 1596. According to Ptolemy there were seven planets, namely Moon, Mercury, Venus, Sun, Mars, Jupiter, Saturn. According to the world-view of Copernicus the Earth and Sun interchanged their roles. The planets were now orbiting around the Sun and not the Earth anymore. Therefore the Earth became a planet, while the Sun was no planet anymore. However, as well the Moon lost its status as a planet since it is orbiting around the Earth. Therefore the number of planets dropped by one from seven to six. How could one explain this number? In his Narratio prima in which he explained Copernicus theory to the world Rheticus gave the following explanation
$$6=1+2+3, \quad 6=1\cdot 2\cdot 3,$$
i.e., $6$ is a perfect number since it equals as well the sum of its divisors as well as its product. Kepler didn't like such an arithmetic explanation. His geometric explanation in the Cosmographic Mystery was the following. If there are six planets there are five intervals between the planets. But five corresponds to the number of Platonic solids. 
\\ \\
To fit the Platonic solids with observational data was challenging, especially after Kepler discovered his first law in his "New Astronomy" (Astronomia Nova)  \cite{kepler1} from 1609 that planets move on ellipses around the sun, which is at rest in a focus of such an ellipse. In his World harmony he was therefore looking for an even deeper explanation for the structure of our world based on consonances. We refer to \cite{field} for a thorough analysis of Kepler's geometrical cosmology.  
\\ \\
According to Kepler there are seven consonant intervals, the octave $\tfrac{1}{2}$, 
the fifth $\tfrac{2}{3}$, the fourth $\tfrac{3}{4}$, the major third $\tfrac{5}{6}$, the minor third $\tfrac{5}{6}$, the major sixth $\tfrac{3}{5}$, and the minor sixth $\tfrac{5}{8}$. That thirds and sixths are considered to be consonant was a rather new trend in Renaissance music 
\cite{walker}. So Kepler distinguished himself from Ptolemy not just that in his cosmology the sun is in the centre instead of the earth, but as well that he considers thirds and sixths as consonant. 
\\ \\
The third book of his Harmonices Mundi is devoted to a geometric explanation why only these seven intervals are consonant. Kepler's idea was that there is a relation between consonance and constructions by ruler and compass. So the 7-gon cannot be constructed by ruler and compass which fits with the fact that $\tfrac{6}{7}$ does not sound consonant. On the other hand, the 15-gon can be constructed by ruler and compass. There Kepler's arguments why the 15-gon should be ruled out look rather ad hoc. Another issue with Kepler's approach is that he was only aware of the constructions using ruler and compass with can be found in Euclid's element. The world still needed to wait for almost two centuries until the 19-year old Gauss discovered how the 17-gon can be constructed by ruler and compass as well.  
\\ \\
In this article based on Kepler's ideas we discuss an approach how Kepler's seven consonant intervals can be characterized mathematically with the help of the numbers of edges of the $n$-gons which were constructed in Euclid's element. We give a mathematical definition which means Euclidean consonant and then prove a theorem that the Euclidean consonants are precisely the Kepler's seven consonant intervals.   
\\ \\
The limited impact of Kepler's geometric approach to the development of music is discussed in the book by Dickreiter \cite{dickreiter}. A particular strong influence Kepler had on Andreas Werckmeister \cite{werckmeister}, the inventor of well tempered tuning, which enabled Johann Sebastian Bach to compose his famous Well-Tempered Clavier.

\section{Euclidean and Gaussian consonances}

On the set of positive rational numbers
$$\mathbb{Q}_+=\big\{q \in \mathbb{Q}: q > 0\big\}$$
we consider the equivalence relation defined for $q_1,q_2 \in \mathbb{Q}_+$ by the requirement
$$q_1 \sim q_2 \quad \Longleftrightarrow \quad q_1=2^n q_2,\,\,n \in \mathbb{Z}.$$
The connection to music is the following. The set $\mathbb{Q}_+$ models rational intervals and two intervals are equivalent if they coincide after octave reduction. The quotient
$$\mathfrak{S}=\mathbb{Q}_+/\sim$$
is than an abelian group with operation induced by the multiplication on $\mathbb{Q}_+$. 
The neutral element in this group
$$\mathfrak{o}=[1]=\bigg[\frac{1}{2}\bigg]$$
than corresponds in the musical interpretation to the octave. 
For $q \in \mathbb{Q}_+$ we denote by $[q] \in \mathfrak{S}$ the equivalence class of $q$.
We refer to the group $\mathfrak{S}$ as the \emph{group of musical intervals} and to its elements as \emph{musical intervals}.
Note that for a musical interval $\sigma \in \mathfrak{S}$ there exist unique positive integers 
$n_\sigma, m_\sigma \in \mathbb{N}$ characterized by the following properties
\begin{equation}\label{rep}
\sigma=\bigg[\frac{m_\sigma}{n_\sigma}\bigg], \quad
\frac{1}{2}\leq\frac{m_\sigma}{n_\sigma} < 1, \quad \mathrm{gcd}(n_\sigma,m_\sigma)=1,
\end{equation}
where $\mathrm{gcd}$ stands for greatest common divisor. In other words, 
the third equality means that $n_\sigma$ and $m_\sigma$ are relative prim. 
\\ \\
We obtain two maps
\begin{eqnarray*}
\mathcal{N} \colon \mathfrak{S} \to \mathbb{N}, \quad \sigma \mapsto n_\sigma\\
\mathcal{M} \colon \mathfrak{S} \to \mathbb{N}, \quad \sigma \mapsto m_\sigma.
\end{eqnarray*}
We further introduce the map
$$K \colon \mathfrak{S} \to \mathfrak{S}, \quad
\sigma \mapsto \bigg[\frac{n_\sigma-m_\sigma}{m_\sigma}\bigg],$$
and refer to it as the \emph{Kepler map}.
We first check that the Kepler map is well-defined. To see this we need to show that
$$\frac{n_\sigma-m_\sigma}{m_\sigma} \in \mathbb{Q}_+.$$
The number is obviously rational. To convince ourselves that it is positive we infer from the second equation in (\ref{rep}) that
$$1<\frac{n_\sigma}{m_\sigma} \leq 2$$
which implies that
\begin{equation}\label{pos}
0<\frac{n_\sigma-m_\sigma}{m_\sigma}\leq 1.
\end{equation}
This implies positivity. Actually for positivity we only need the first inequality in (\ref{pos}). 
\\ \\
The following lemma tells us that the map $\mathcal{N}$ is monotone decreasing under the Kepler map. 
\begin{lemma}\label{mono}
For every $\sigma \in \mathfrak{S}$ we have $\mathcal{N} \circ K(\sigma) \leq \mathcal{N}(\sigma)$
and equality holds iff $\sigma=\mathfrak{o}$. 
\end{lemma}
\textbf{Proof: } We first consider the case $\sigma \neq \mathfrak{o}$. In this case we have
$$\frac{1}{2}<\frac{m_\sigma}{n_\sigma}$$
implying that
$$n_\sigma<2m_\sigma$$
and therefore
$$n_\sigma-m_\sigma<m_\sigma$$
so that 
$$\frac{n_\sigma-m_\sigma}{m_\sigma}<1.$$
Hence there exists a unique nonnegative ineger $\ell \in \mathbb{N}_0$ such that
\begin{equation}\label{k1}
\frac{1}{2} \leq \frac{2^\ell(n_\sigma-m_\sigma)}{m_\sigma}<1.
\end{equation}
Decompose $m_\sigma$ uniquely into the product of an odd number and a power of $2$, i.e.,
$$m_\sigma=2^k \mu_\sigma$$
where $k \in \mathbb{N}_0$ and $\mu_\sigma$ is odd. 
Define
$$\rho:=\max\{0,k-\ell\} \in \mathbb{N}_0, \qquad \tau:=\max\{0,\ell-k\} \in \mathbb{N}_0.$$
We claim 
\begin{equation}\label{kf}
n_{K(\sigma)}=2^\rho \mu_\sigma, \qquad m_{K(\sigma)}=2^{\tau}(n_\sigma-m_\sigma).
\end{equation}
To prove (\ref{kf}) we need according to (\ref{rep}) check the following
\begin{eqnarray}\label{kf1}
K(\sigma)=\bigg[\frac{2^{\tau}(n_\sigma-m_\sigma)}{2^\rho \mu_\sigma}\bigg],\\ \nonumber
\frac{1}{2} \leq \frac{2^{\tau}(n_\sigma-m_\sigma)}{2^\rho \mu_\sigma}<1,\\ \nonumber
\mathrm{gcd}\big(2^\rho \mu_\sigma, 2^\tau(n_\sigma-m_\sigma)\big)=1.
\end{eqnarray}
By shortening the fraction we obtain
$$\frac{2^\ell(n_\sigma-m_\sigma)}{m_\sigma}=\frac{2^\ell(n_\sigma-m_\sigma)}{2^k \mu_\sigma}=\frac{2^{\tau}(n_\sigma-m_\sigma)}{2^\rho \mu_\sigma}.$$
Hence the second assertion in (\ref{kf1}) follows from (\ref{k1}). Moreover, the first assertion in (\ref{kf1}) follows from
$$K(\sigma)=\bigg[\frac{n_\sigma-m_\sigma}{m_\sigma}\bigg]=
\bigg[\frac{2^\ell(n_\sigma-m_\sigma)}{m_\sigma}\bigg]=\bigg[
\frac{2^{\tau}(n_\sigma-m_\sigma)}{2^\rho \mu_\sigma}\bigg].$$
It remains to check the third assertion in (\ref{kf1}). First note that since 
$n_\sigma$ and $m_\sigma$ are relative prim, the same is true for $m_\sigma$ and
$n_\sigma-m_\sigma$ so that we have
\begin{equation}\label{relprim0}
\mathrm{gcd}(m_\sigma,n_\sigma-m_\sigma)=1.
\end{equation}
Since $\mu_\sigma$ is a divisor of $m_\sigma$ this implies
\begin{equation}\label{relprim}
\mathrm{gcd}(\mu_\sigma,n_\sigma-m_\sigma)=1.
\end{equation}
We now consider first the case where $\rho=0$. Since $\mu_\sigma$ is odd the third assertion in (\ref{kf1}) follows in this case directly from (\ref{relprim}). It remains to check the case where 
$$\rho\neq 0.$$
This implies that
$$\tau=0, \quad k>0.$$
In particular, since $k>0$ we conclude that $m_\sigma=2^k \mu_\sigma$ is even. From
(\ref{relprim0}) we conclude that $n_\sigma-m_\sigma$ is odd. Hence again the third assertion in (\ref{kf1}) follows from (\ref{relprim}). Hence the truth of all three assertion in (\ref{k1}) is establishes and (\ref{kf}) follows. 
\\ \\
From the first assertion in (\ref{kf}) we conclude 
\begin{equation}\label{kepseq}
\mathcal{N} \circ K(\sigma)=n_{K(\sigma)}=2^\rho \mu_\sigma \leq 2^k \mu_\sigma=m_\sigma<n_\sigma=\mathcal{N}(\sigma)
\end{equation}
where the last inequality follows from the second property in (\ref{rep}). This proves the lemma in the case where $\sigma \neq \mathfrak{o}$.
\\ \\
It remains to discuss the case where $\sigma=\mathfrak{o}$. We have
$$n_\mathfrak{o}=2, \quad m_\mathfrak{o}=1.$$
It follows that
\begin{equation}\label{fix}
K(\mathfrak{o})=\bigg[\frac{n_\mathfrak{o}-m_\mathfrak{o}}{m_\mathfrak{o}}\bigg]=[1]=\mathfrak{o}
\end{equation}
so that $\mathfrak{o}$ is a fixed point of $K$. In particular,
$$\mathcal{N} \circ K(\mathfrak{o})=\mathcal{N}(\mathfrak{o}).$$
This finishes the proof of the lemma. \hfill $\square$
\begin{cor}\label{monocor}
Given $\sigma \in \mathfrak{S}$ there exists $\ell \in \mathbb{N}$ such that
$K^\ell(\sigma)=\mathfrak{o}$. 
\end{cor}
\textbf{Proof: } As we have seen in (\ref{fix}) the octave is a fixed point of the Kepler map so that we can assume without loss of generality that $\sigma \neq \mathfrak{o}$. By Lemma~\ref{mono} there
exists $\ell \in \mathbb{N}$ so that 
$$\mathcal{N} \circ K^\ell(\sigma) \leq 2.$$
By the second assertion in (\ref{rep}) it follows that 
$$\mathcal{N} \circ K^\ell(\sigma) = 2, \quad \mathcal{M} \circ K^\ell(\sigma)=1$$
so that 
$$K^\ell(\sigma)=\mathfrak{o}.$$
This proves the Corollary. \hfill $\square$ 
\\ \\
In view of Corollary~\ref{monocor} we can make the following definition.
\begin{fed}
We define the \emph{height} of a musical interval $\sigma \in \mathfrak{S}$ as
$$\mathfrak{h}(\sigma):=\min\big\{\ell \in \mathbb{N}_0: K^\ell(\sigma)=\mathfrak{o}\big\}.$$
\end{fed}
Note that using the definition of height, the octave $\mathfrak{o}$ can be characterized as the unique musical interval of height $0$.
\\ \\
Given a musical interval $\sigma \in \mathfrak{S}$ we define its \emph{first Kepler sequence} as
$$\mathfrak{K}_1(\sigma):=\big(\sigma,K(\sigma),\ldots,K^{\mathfrak{h}(\sigma)}(\sigma)\big\},
$$
and its \emph{second Kepler sequence} as
$$\mathfrak{K}_2(\sigma):=\mathcal{N}(\mathfrak{K}_1(\sigma)):=
\big(\mathcal{N}(\sigma),\mathcal{N} \circ K(\sigma),\ldots, \mathcal{N} \circ
K^{\mathfrak{h}(\sigma)}(\sigma)\big).$$
We consider some examples. The first Kepler sequence for the minor sixth $\big[\tfrac{5}{8}\big]$ is
\begin{equation}\label{k1sext}
\mathfrak{K}_1\Bigg(\bigg[\frac{5}{8}\bigg]\Bigg)=
\Bigg(\bigg[\frac{5}{8}\bigg],\bigg[\frac{3}{5}\bigg],\bigg[\frac{2}{3}\bigg],
\bigg[\frac{1}{2}\bigg]\Bigg),
\end{equation}
i.e.
$$\big(\textrm{minor sixth}, \textrm{major sixth}, \textrm{fifth}, \textrm{octave}\big).$$
The second Kepler sequence for the minor sixth is than
\begin{equation}\label{k2sext}
\mathfrak{K}_2\Bigg(\bigg[\frac{5}{8}\bigg]\Bigg)=
\big(8,5,3,2\big).
\end{equation}
For the minor third $\big[\tfrac{5}{6}\big]$ we get as first Kepler sequence
\begin{equation}\label{k1terz}
\mathfrak{K}_1\Bigg(\bigg[\frac{5}{6}\bigg]\Bigg)=\Bigg(\bigg[\frac{5}{6}\bigg],
\bigg[\frac{4}{5}\bigg],\bigg[\frac{1}{2}\bigg]\Bigg),
\end{equation}
i.e.,
$$\big(\textrm{minor third},\textrm{major third},\textrm{octave}\big),$$
and as second Kepler sequence we have
\begin{equation}\label{k2terz}
\mathfrak{K}_2\Bigg(\bigg[\frac{5}{6}\bigg]\Bigg)=\big(6,5,2\big).
\end{equation}
For the fourth $\big[\tfrac{3}{4}\big]$ the first Kepler sequence is
\begin{equation}\label{k1quarte}
\mathfrak{K}_1\Bigg(\bigg[\frac{3}{4}\bigg]\Bigg)=\Bigg(\bigg[\frac{3}{4}\bigg],
\bigg[\frac{2}{3}\bigg],\bigg[\frac{1}{2}\bigg]\Bigg),
\end{equation}
i.e.,
$$\big(\textrm{fourth},\textrm{fifth},\textrm{octave}\big)$$
and the second Kepler sequence
\begin{equation}\label{k2quarte}
\mathfrak{K}_2\Bigg(\bigg[\frac{3}{4}\bigg]\Bigg)=\big(4,3,2\big).
\end{equation}
Inspired by Kepler our next goal is to define a musical interval as consonant if the members of its second Kepler sequence are integers which give rise to constructible polygons. For this purpose recall that a prime number $p$ is called \emph{Fermat prime} if it is of the form
$$p=2^{2^n}+1.$$
There are five known Fermat primes
\begin{equation}\label{fermat}
F_0=3,\quad F_1=5, \quad F_2=17, \quad F_3=257, \quad F_4=65537
\end{equation}
and there are some probabilistic considerations that these might be the only
ones \cite{boklan-conway}. The Gauss-Wantzel theorem tells us that a regular $n$-gon is constructible by ruler and compass if and only if $n$ is the product of a power of $2$ and any number of distinct Fermat primes \cite{gauss,wantzel}. Gauss proved that these $n$-gons can be constructed while Wantzel showed that no other $n$-gon is constructible. We make the following definition.
\begin{fed}
A positive integer $n>1$ is called a \emph{Gauss-Wantzel number} if
$$n=2^k p_1\cdots p_\ell$$
where $k,\ell \in \mathbb{N}_0$ and $p_1,\ldots,p_\ell$ are distinct Fermat primes.  
\end{fed}
We introduce 
$$\mathfrak{G}:=\big\{n \in \mathbb{N}: n\,\,\textrm{Gauss-Wantzel}\big\}$$
the subset of Gauss-Wantzel numbers. 
With this notion the Gauss-Wantzel theorem just tells us that an $n$-gon is constructible with ruler and compass iff $n \in \mathfrak{G}$. 
\\ \\
In Euclid's Elements \cite{euclid} one can find the construction of $n$-gons using ruler and compass in the special case where the only allowed Fermat primes are $3$ and $5$. Hence we make the following definition.
\begin{fed}
A Gauss-Wantzel number is called \emph{Euclidean} if 
$$n=2^k 3^\ell 5^m$$
where $k \in \mathbb{N}_0$ and $\ell,m \in \{0,1\}$. 
\end{fed}
The first eleven Euclidean Gauss-Wantzel numbers are
$$2,3,4,5,6,8,10,12,15,16,20,\ldots$$
The smallest Gauss-Wantzel number which is not Euclidean is 17. It was a big surprise when the 19 year old Gauss announced 1796 in the ``Intelligenzblatt der allgemeinen Literaturzeitung" a construction of the 17-gon by ruler and compass \cite{gauss0}. Five years later he published his construction in his ``Disquisitiones Arithmeticae" \cite{gauss}.
\\ \\
We abbreviate
$$\mathfrak{E}:=\big\{n \in \mathfrak{G}: n\,\,\textrm{Euclidean}\big\}$$
the set of Euclidean Gauss-Wantzel numbers. We therefore have a nested sequence of subsets
$$\mathfrak{E} \subset \mathfrak{G} \subset \mathbb{N}.$$
We are now in position to define two versions of consonant musical intervals.
\begin{fed}
A musical interval $\sigma \in \mathfrak{S}$ is called \emph{Euclidean consonant} if
$$\mathfrak{K}(\sigma) \in \mathfrak{E}^{\mathfrak{h}(\sigma)+1},$$
i.e., all members of its second Kepler sequence are Euclidean Gauss-Wantzel numbers. 
\end{fed}
\begin{fed}
A musical interval $\sigma \in \mathfrak{S}$ is called \emph{Gaussian consonant} if
$$\mathfrak{K}(\sigma) \in \mathfrak{G}^{\mathfrak{h}(\sigma)+1},$$
i.e., all members of its second Kepler sequence are Gauss-Wantzel numbers. 
\end{fed}
Since $\mathfrak{E} \subset \mathfrak{G}$ every Gaussian consonant is as well a Euclidean consonant. Further note, that if $\sigma$ is an Euclidean respectively Gaussian consonant any musical interval in the first Kepler sequence $\mathfrak{K}_1(\sigma)$ is as well an Euclidean
respectively Gaussian consonant. 
\\ \\
We have now the following theorem which goes back to the ideas of Kepler.
\begin{thm}
There are precisely seven Euclidean consonants, namely the minor third, the major third, the fourth, the fifth, the minor sixth, the major sixth, and the octave, i.e, 
$$\Bigg\{\bigg[\frac{5}{6}\bigg], \bigg[\frac{4}{5}\bigg],\bigg[\frac{3}{4}\bigg],
\bigg[\frac{2}{3}\bigg], \bigg[\frac{5}{8}\bigg],\bigg[\frac{3}{5}\bigg],\bigg[\frac{1}{2}\bigg]\Bigg\} \subset
\mathfrak{S}.$$
\end{thm}
\textbf{Proof: } That this seven musical intervals are Euclidean consonants follows from (\ref{k1sext})--(\ref{k2quarte}). It remains to show that there are no others. Suppose that
$$\sigma=\bigg[\frac{m_\sigma}{n_\sigma}\bigg]$$
is a Euclidean consonant. We necessarily have
$$n_\sigma \in \mathfrak{E}.$$ 
From (\ref{kepseq}) in the prove of Lemma~\ref{mono} we conclude that there exists 
$\ell \in \mathbb{N}_0$ such that
$$2^\ell n_{K(\sigma)}=m_\sigma.$$
Since $\sigma$ is an Euclidean consonant, it holds that
$$n_{K(\sigma)} \in \mathfrak{E}.$$
Since $\mathfrak{E}$ is invariant under multiplication by $2$ we conclude that
$$2^\ell n_{K(\sigma)} \in \mathfrak{E}$$ and therefore
$$m_\sigma \in \mathfrak{E}.$$
We first discuss the case where $n_\sigma$ is odd, i.e., 
$$n_\sigma \in \{3,5,15\}.$$ 
By the second assertion in (\ref{rep}) we have 
\begin{equation}\label{ineq}
\frac{n_\sigma}{2}<m_\sigma<n_\sigma.
\end{equation}
Hence for $n_\sigma=3$ the only possibility is the fifth $\big[\tfrac{2}{3}\big]$ and for
$n_\sigma=5$ the only possibilities are the major sixth $\big[\tfrac{3}{5}\big]$ and the major third $\big[\tfrac{4}{5}\big]$. In the odd case it therefore suffices to discuss the case 
$n_\sigma=15$. Since $n_\sigma$ and $m_\sigma$ are relatively prime by the third assertion in
(\ref{rep}) we do not need to worry about the cases $m_\sigma \in \{9,10,12\}$. So the only case left to discuss in the odd case is
$$\sigma=\bigg[\frac{8}{15}\bigg].$$
We compute
$$K\Bigg(\bigg[\frac{8}{15}\bigg]\Bigg)=\bigg[\frac{7}{8}\bigg],\quad
K^2\Bigg(\bigg[\frac{8}{15}\bigg]\Bigg)=\bigg[\frac{4}{7}\bigg]$$
so that
\begin{equation}\label{nichtschoen}
\mathcal{N} \circ K^2\Bigg(\bigg[\frac{8}{15}\bigg]\Bigg)=7 \notin \mathfrak{E}.
\end{equation}
We conclude that $\big[\tfrac{8}{15}\big]$ is not Euclidean consonant. In particular, there are no further Euclidean consonants in the odd case.
\\ \\
It remains to discuss the case where $n_\sigma$ is even, i.e, $n_\sigma=2^k3^\ell 5^m$ for $k \in \mathbb{N}$ and $m,k \in \{0,1\}$. Since $n_\sigma$ and $m_\sigma$ are relativ prime it follows that $m_\sigma$ has to be odd, i.e., $m_\sigma \in \{3,5,15\}$. In view of (\ref{ineq}) the only case which appears for $m_\sigma=3$ is the fourth $\big[\tfrac{3}{4}\big]$. 
For $m_\sigma=5$ there are only the minor third $\big[\tfrac{5}{6}\big]$ and the minor
sixth $\big[\tfrac{5}{8}\big]$. If $m_\sigma=15$, then since $n_\sigma$ and $m_\sigma$ are relativ prime it follows that $n_\sigma=2^k$ is a power of $2$. In view of (\ref{ineq}) the only case left to discuss is $n_\sigma=16$, i.e., the musical interval $\big[\tfrac{15}{16}\big]$. We compute
$$K\Bigg(\bigg[\frac{15}{16}\bigg]\Bigg)=\bigg[\frac{8}{15}\bigg].$$
We have already seen in (\ref{nichtschoen}) that $\big[\tfrac{8}{15}\big]$ is not Euclidean consonant. Therefore $\big[\tfrac{15}{16}\big]$ is not Euclidean consonant either. This finishes the proof of the theorem. \hfill $\square$
\\ \\
Apart from the Euclidean consonants there are additional Gaussian consonants. We have the following lemma.
\begin{lemma}
If $F$ is a Fermat prime, then $\big[\frac{F-1}{F}\big]$ is a Gaussian consonant.
\end{lemma}
\textbf{Proof: } A Fermat prime is a prime number of the form $F=2^{2^n}+1$ for $n \in \mathbb{N}_0$. Therefore we have
$$\frac{F-1}{F}=\frac{2^{2^n}}{2^{2^n}+1},$$
so that we get
$$K\Bigg(\bigg[\frac{F-1}{F}\bigg]\Bigg)=\bigg[\frac{1}{2^{2^n}}\bigg]=\bigg[
\frac{1}{2}\bigg].$$
Therefore we obtain as first Kepler sequence
$$\mathfrak{K}_1\Bigg(\bigg[\frac{F-1}{F}\bigg]\Bigg)=\Bigg(\bigg[\frac{F-1}{F}\bigg],
\bigg[\frac{1}{2}\bigg]\Bigg),$$
and as second Kepler sequence
$$\mathfrak{K}_2\Bigg(\bigg[\frac{F-1}{F}\bigg]\Bigg)=\big(F,2\big),$$
all whose members are Gauss-Wantzel numbers. This proves the lemma. \hfill $\square$
\\ \\
In view of the known Fermat primes (\ref{fermat}) we obtain apart from the fifth
$\big[\tfrac{2}{3}\big]$ and the major third $\big[\tfrac{4}{5}\big]$ as Gaussian consonants
$$\bigg[\frac{F_2-1}{F_2}\bigg]=\bigg[\frac{16}{17}\bigg], \quad
\bigg[\frac{F_3}-1{F_3}\bigg]=\bigg[\frac{256}{257}\bigg], \quad
\bigg[\frac{F_4-1}{F_4}\bigg]=\bigg[\frac{65536}{65537}\bigg].$$
These three Gaussian consonants are not Euclidean consonants. Another example of a Gaussian consonant which is not Euclidean consonant is
$$\sigma=\bigg[\frac{12}{17}\bigg].$$
Its image under the Kepler map
$$K(\sigma)=\bigg[\frac{5}{6}\bigg]$$
is the minor third, which is Euclidean consonant. Therefore since
$$\mathcal{N}(\sigma)=17 \in \mathfrak{G}$$
it follows that $\sigma$ is Gaussian consonant.


\begin{thebibliography}{99}

\bibitem{boklan-conway}
K.\,Boklan, J.\,Conway \emph{Expect at most one billionth of a new Fermat prime!},
Math.\,Intelligencer \textbf{39}, no.\,1 (2017), 3--5.

\bibitem{dickreiter} M.\,Dickreiter, \emph{Der Musiktheoretiker Johannes Kepler},
Bern; M\"unchen, Francke (1973)

\bibitem{euclid} Euclid, Euclid's Elements, Sir Thomas Little Heath, New York, Dover (1956)

\bibitem{field} J.\,Field, \emph{Kepler's geometrical cosmology}, Chicago, Univ.\,of Chicago Pr. (1987).

\bibitem{gauss0} C.\,Gauss, \emph{Neue Entdeckungen}, Intelligenzblatt der allgem.\,Literatur-Zeitung \textbf{66} (1796), 544.

\bibitem{gauss} C.\,Gauss, \emph{Disquisitiones Arithmeticae}, Leipzig, Gerh.\,Fleischer (1801).

\bibitem{kepler0} J.\,Kepler translated by E.\,Aiton, A.\,Duncan \emph{Mysterium Cosmographicum. The Secret of the Universe}, The second edition of Kepler's work, reprinted with translation by A.\,Duncan and introduction and commentary by E.\,Aiton, New York (1981).

\bibitem{kepler1} J.\,Kepler translated by W.\,Donahue \emph{New Astronomy}, Cambridge: Cambridge Univ.\,Pr. (1992).

\bibitem{kepler2} J.\,Kepler translated to German by Max Caspar, \emph{Weltharmonik}, M\"unchen-Berlin (1939).

\bibitem{walker}
D.\,Walker, \emph{Kepler's Celestial Music}, J.\,of the Warburg and Courtauld Institutes \textbf{30} (1967), 228--250.

\bibitem{wantzel}
L.\,Wantzel, \emph{Recherches sur les moyens de reconna\^itre si un Probl\`eme de G\'eometrie
peut se r\'esoudre avec la r\`egle et le compas}, J.\,de Math.\,Pures et Appliqu\'ees
\textbf{2} (1837), 366--372.

\bibitem{werckmeister}
A.\,Werckmeister, \emph{Musicalische Paradoxal-Discourse}, Quedlinburg (1707). 


\end{thebibliography}
\end{document}